\documentclass[letterpaper,11pt]{article}
\usepackage{amsfonts,latexsym,amssymb,amsmath,algorithmic,epsfig,rotating}

\newtheorem{thm}{Theorem}

\newtheorem{lemma}{Lemma} 
\newtheorem{rem}{Remark}
\newtheorem{dfn}{Definition}

\newtheorem{prop}{Proposition}

\newlength{\jmr}
\newlength{\khov}
\newlength{\bernd}
\newcommand{\thth}{$^{\text{\underline{th}}}${} } 
\newcommand{\bla}{{\overline{\nabla}_\cA}}

\newcommand{\Pro}{{\mathbb{P}}}

\newcommand{\R}{\mathbb{R}}
\newcommand{\C}{\mathbb{C}}

\newcommand{\Z}{\mathbb{Z}}
\newcommand{\bO}{\mathbf{O}}
\newcommand{\cA}{\mathcal{A}}

\newcommand{\Zn}{\Z^n}

\newcommand{\Rn}{\R^n}

\newcommand{\Cs}{\C^*}
\newcommand{\Rs}{\R^*}

\newcommand{\cT}{{\mathcal{T}}}

\newcommand{\Rsn}{{(\R^*)}^n}
\newcommand{\Csn}{{(\C^*)}^n}

\newcommand{\qed}{$\blacksquare$}
\newcommand{\dia}{$\diamond$}

\newcommand{\supp}{\mathrm{Supp}}
\newcommand{\conv}{\mathrm{Conv}}

\newcommand{\codim}{\mathrm{codim }}

\begin{document}
\title{\mbox{}\\
\vspace{-1in}
New Complexity Bounds for Certain Real Fewnomial Zero Sets\\ 
(Extended Abstract)} 

\author{
Joel Gomez\thanks{ Department of Mathematics, Box 8205, NC State University, 
Raleigh, NC 27695-8205. Partially supported by an undergraduate fellowship 
from NSF REU grant DMS-0552610. } \and 
Andrew Niles\thanks{ CPU Box 274436, University of Rochester, 
Rochester, NY \ 14627.  Partially supported 
by an undergraduate fellowship from NSF REU grant DMS-0552610. } \and  
J.\ Maurice Rojas\thanks{
Department of Mathematics,  
Texas A\&M University
TAMU 3368, 
College Station, Texas \ 77843-3368,  
USA,  {\tt rojas@math.tamu.edu} \ ,  
{\tt www.math.tamu.edu/\~{}rojas} \ .
Partially supported by NSF CAREER grant DMS-0349309 and 
NSF REU grant DMS-0552610.}   } 

\date{\today} 

\maketitle

\mbox{}\hfill
{\em Rojas dedicates this paper to his friend, Professor Tien-Yien Li. } 
\hfill\mbox{}

\begin{abstract} 
Consider real bivariate polynomials $f$ and $g$, respectively 
having $3$ and $m$ monomial terms. We prove 
that for all $m\!\geq\!3$, there are systems of the form $(f,g)$ having 
exactly $2m-1$ roots in the positive quadrant. 
Even examples with $m\!=\!4$ having $7$ positive roots 
were unknown before this paper, so we detail 
an explicit example of this form. We also present an $O(n^{11})$ 
upper bound for the number of diffeotopy types of the real zero 
set of an $n$-variate polynomial with $n+4$ monomial terms. 
\end{abstract} 

\section{Introduction} 
Finding the correct combinatorics governing the real zero sets  
of sparse polynomials is a major open problem within real algebraic geometry 
(see, e.g., \cite{bs,drrs}). In particular, while the maximal number of real 
roots of a sparse polynomial 
in one variable has been well-understood for centuries (dating back to 
17\thth{} century work of Descartes \cite{descartes}), only loose 
upper bounds are known in higher-dimensions. Nevertheless, the bounds 
currently known have already proved of great use in arithmetic 
geometry \cite{cohenzannier} and Hilbert's 16\thth{} Problem \cite{kaloshin} 
(to name just a few areas), 
and it is known that {\bf optimal} bounds would have significant applications 
in many areas of engineering. 
Here, with an eye toward tightening known upper bounds, we exhibit 
sparse polynomial systems with more roots than previously known 
(Theorem \ref{thm:and} below), 
and a new {\bf polynomial} upper bound on the number of diffeotopy types of 
real zero sets of $n$-variate polynomials with $n+4$ monomial terms 
(Theorem \ref{thm:3d} below). 

\subsection{New Lower Bounds} 
Consider the following system of analytic equations: 
\[ (\star) \ \ \left\{ \begin{matrix} 
\alpha_1+\alpha_2x^{a_2}y^{b_2}+\alpha_3x^{a_3}y^{b_3}\\
\beta_1 + \beta_2 x^{c_2}y^{d_2} + \cdots +\beta_m x^{c_m} y^{d_m},
\end{matrix}\right. \]   
for nonzero real $\alpha_i$ and $\beta_i$ and distinct nonzero real vectors 
$(a_2,b_2)$, $(a_3,b_3)$, $(c_2,d_2),\ldots,(c_m,d_m)$. 
The number of isolated roots $(x,y)$ in the positive quadrant 
$\R^2_+$ is of course 
bounded above by some function of the coefficients and exponents, but it 
wasn't until Askold Khovanski's invention of {\bf Fewnomial Theory} 
around the early 1980s \cite{kho,few} that an upper bound depending solely 
on $m$ was known: $3^{m+2}2^{(m+1)(m+2)/2}$ (just for the number of 
non-degenerate roots),\footnote{Via now standard tricks, Khovanski's 
bound easily implies an upper bound of $2^{O(m^2)}$ for the number 
of isolated roots.} invoking a very special 
case of his famous {\bf Theorem on Real Fewnomials}. 

Later, Li, Rojas, and Wang proved an upper bound of $2^m-2$ 
\cite{lrw}, while more recently Avenda\~{n}o has proved that 
for the special case where $(a_2,b_2,a_3,b_3)\!=\!(1,0,0,1)$ --- and the 
remaining $(c_i,d_i)$ lie in $\Zn$ --- the {\em polynomial} system ($\star$) 
never has more than $2m-2$ isolated positive roots \cite{avendano}. It has been 
conjectured that the correct general upper 
bound for the number of isolated positive roots of ($\star$) should be 
polynomial in $m$, but this remains an open 
problem. So we provide the following new lower bound. 
\begin{thm}
\label{thm:and} 
For all $m\!\geq\!3$, there exist polynomial systems of the form ($\star$) 
above with at least $2m-1$ roots in $\R^2_+$. In particular, the polynomial 
system 
\[ \begin{matrix}x^6+\frac{44}{31}y^3-y\\ 
y^{14}+\frac{44}{31}x^3y^8-xy^8+\alpha x^{133}
\end{matrix} \] 
has exactly $7$ roots in $\R^2_+$ for $1936838\!\leq\!\alpha\!\leq\!1936254$. 
\end{thm} 

\noindent 
While it easy to construct systems of the form ($\star$) with exactly $2m-2$ 
positive roots, examples with $m\!=\!3$ having $5$ positive roots weren't known
until 2000 \cite{haas,drrs}. Moreover, even examples with $m\!=\!4$ 
having $7$ positive roots appear to have been unknown before this paper. 

\subsection{New Topological Upper Bounds} 
Recall that while a smooth, real, degree $d$ projective plane curve has
at most $1+\text{\scalebox{.7}[.7]{$\begin{pmatrix}d-1\\ 2\end{pmatrix}$}}$
connected components \cite{harnack}, determining the possible nestings of
these ovals --- a piece of the first part of Hilbert's famous 16\thth{} 
Problem \cite{kaloshin} --- is quite complicated. In more general language,
this is the determination of possible {\bf diffeotopy types} of
such curves.
\begin{dfn}
Recall that a {\bf diffeotopy} between two sets
$X,Y\!\subseteq\!\Rn$ is a differentiable function $H : [0,1]\times \Rn
\longrightarrow \Rn$
such that $H(t,\cdot)$ is a diffeomorphism for all $t\!\in\![0,1]$,
$H(0,\cdot)$ is the identity on $X$, and $H(1,X)\!=\!Y$. Equivalently,
we simply say that $X$ and $Y$ are {\bf diffeotopic}. \dia
\end{dfn}
                                                                                
\noindent
Note that diffeotopy is a more refined equivalence than diffeomorphism, since
diffeotopy implies an entire continuous family of ``infinitesimal''
diffeomorphisms that deform $X$ to $Y$ and back again. Returning to
nestings of ovals of real degree $d$ projective plane curves, an asymptotic
formula of $e^{d^2}$ is now known \cite{orevkov}, and the {\bf exact} number is
currently known (as of early 2007) only for $d\!\leq\!8$.
                                                                                
Via our techniques here, we can count diffeotopy types in a dramatically
different setting. Recall that $\R^*\!:=\!\R\!\setminus\{0\}$. 
\begin{dfn}
Given any polynomial $f$, its {\bf support} (or
{\bf spectrum}) --- written $\supp(f)$ --- is the set of exponent vectors 
in its monomial terms.
Also, we let $Z_+(f)$ (resp.\ $Z^*_\R(f)$) denote the set of\linebreak 
\scalebox{.97}[1]{roots of $f$ in
$\Rn_+$ (resp.\ $\Rsn$). Finally, given any $\cA\!\subset\!\Rn$, we let 
$\conv \cA$ denote its {\bf convex hull}. \dia} 
\end{dfn}
\begin{thm}
\label{thm:3d}
For any fixed $\cA\!\subset\!\Zn$
with $\#\cA\!=\!n+4$ and $\conv \cA$ of positive volume, 
there are no more than $O(n^{11})$
diffeotopy types for any smooth
$Z^*_\R(f)$ with $\supp(f)\!=\!\cA$. In particular, our bound 
is completely independent of the coordinates of $\cA$.
\end{thm}

\noindent
The positive volume assumption is natural, for otherwise one would in fact be
studying (up to an invertible monomial change of variables) an
instance where $\#\cA\!\geq\!n'+5$ with $n'\!\leq\!n$ and $\conv \cA$ 
has positive $n'$-dimensional volume (see, e.g., \cite[Cor.\ 1,  
Sec.\ 2]{thresh}).
                                                                 
That there exists any sort of upper bound depending solely on $n$ is already a 
non-trivial fact, first observed by
Lou van den Dries around the 1990s via o-minimality (see, e.g.,
\cite[Prop.\ 3.2, Pg.\ 150]{ominimal}). Bounds exponential in 
$n$ were then (implicitly) discovered in \cite{gabrielov}, and 
appear explicitly in \cite{basu}. 
Our polynomial bound above is thus a great improvement. 
In particular, {\bf polynomial} bounds for the number of diffeotopy types 
were previously known only for $\#\cA\!\leq\!n+3$ \cite{thresh,drrs}. 

\section{Outlines of the Proofs} 
The complete proofs of our two main theorems will appear in the 
full version of this paper. In this extended abstract, we will simply sketch 
the main ideas. 
\subsection{Proving Theorem \ref{thm:and} by Induction} 
The special case $m\!=\!3$ follows immediately from earlier 
work of Haas \cite{haas}, where an example consisting of 
two degree $106$ polynomials was detailed. The more recent 
paper \cite{drrs} gives a far simpler example, consisting of a 
pair of degree $6$ polynomials, and gives an explanation of 
the paucity of such extremal examples via $\cA$-discriminants. 

The special case $m\!=\!4$ can be checked by computationally 
verifying our stated example, e.g., via rational univariate reduction and 
an application of Sturm-Habicht 
sequences to count real roots. This approach is dates back to the 
19\thth century and has undergone recent algorithmic revivals 
\cite{kron,pspace,gls}, so the 
verification of our example is a simple exercise in {\tt Maple}   
(see \cite{drrs} for an extended illustration of these and other techniques 
when $m\!=\!3$). The general construction of such examples is more subtle, 
however, so let us now assume $m\!\geq\!5$. 

By rescaling the variables, dividing by suitable monomial terms, 
and employing a monomial change of variables \cite{lrw}, we 
can reduce to the following univariate problem: Find real 
$a_i,b_j,c_k$ such that 
\[ 1+c_2x^{a_2}(1-x)^{b_2}+\cdots+c_m x^{a_m}(1-x)^{b_m}\] 
has at least $2m-1$ roots in the open interval $(0,1)$. 

We now proceed by induction on $m$: Let us 
assume that we have constructed our desired 
example for some fixed $m\!\geq\!4$, and that it in fact 
has the following form: 
\[ f_{m}(x):=1-\left(\frac{31}{44}\right)^{35/12}x^{-1/6}(1-x)^{35/12}-
\left(\frac{44}{31}\right)^{5/6}x^{1/3}(1-x)^{1/6}-
c_4x^{a_4}(1-x)^{b_4}-\cdots - c_{m}x^{a_{m}}(1-x)^{b_{m}}, \] 
for some positive integers $a_4,b_4,\ldots,a_m,b_m$. 
We will then show that we can find a positive real $c$ and a 
positive integer $a$ such that 
\[ g(x):=f_m(x)-cx^{a}(1-x)^{7}\] has at least $2m+1$ roots in $(0,1)$. 

Toward this end, first note that (a) for any $m$, the function $f_m$ satisfies 
$\lim\limits_{x\rightarrow 0^+} f_m(x)\!=\!-\infty$ and 
$f_m(1)\!=\!1$, and (b) the graph of $x^{a}(1-x)^7$ is an arc in 
the first quadrant, {\bf nearly flat aside from a single peak}, connecting 
$(0,0)$ and $(1,0)$. In particular, the graph of $x^{a}(1-x)^7$ 
becomes ``flatter'' as $a$ increases. The main trick will then be 
to pick $a$ sufficiently large (and $c\!>\!0$ a bit more carefully) so that 
the graphs of $f_m$ and $x^{a_{m+1}}(1-x)^7$ intersect in at 
least $2m+1$ points in the first quadrant.  

Labelling the roots of $f_m$ in $(0,1)$ as $z_1\!<\cdots<\!z_{2m-1}$, 
we then need only make some elementary calculations 
to find conditions on $c$ and $a$ so that guarantee that our graph condition 
holds. Toward this end, let $y_1$ (resp.\ $y_{2m}$ and $y_{2m+1}$) be the 
unique 
point in $(0,z_1)$ (resp.\ $(z_{2m-1},1)$ and $(y_{2m},1)$) such 
that $f_m(y_1)$ is $-1/4$ (resp.\ $1/4$ and $1/2$). 
Also, by Rolle's Theorem, for all $i\!\in\!\{2,\ldots,2m-1\}$, we can find a 
point $y_i\!\in\!(z_{i-1},z_i)$ such that $f'_m(y_i)\!=\!0$. 

In order to enforce our graph condition, it is easily checked (employing 
Rolle's Theorem once more) that 
it suffices to find conditions on $c$ and $a$ that imply 
$cy^a_{2m}(1-y_{2m})^7\!<\!\min\limits_{1<i<2m}|f_m(y_i)|$ and $cy^a_{2m+1}
(1-y_{2m+1})^7\!>\!f_m(y_{2m+1})\!=\!1/2$. Letting 
$\alpha\!:=\!\min\limits_{1<i<2m}|f_m(y_i)|$, an elementary calculation  
then yields the sufficients conditions  
$cy^a_{2m}(1-y_{2m})^7\!<\!\alpha$ and $cy^a_{2m+1}(1-y_{2m+1})^7\!>\!1/2$. 

With a bit more work, we then obtain that any $c$ satisfying 
$\frac{\alpha}{y^a_{2m}(1-y_{2m})^7}\!>\!c\!>\!\frac{1}{2y^a_{2m+1}
(1-y_{2m+1})^7}$ will work, provided that 
$a\!>\!\left.\log\left(\frac{1}{2\alpha}\left(\frac{1-y_{2m}}{1-y_{2m+1}}
\right)^7\right)\right/\log\left(\frac{y_{2m+1}}{y_{2m}}\right)$. 
Since $y_{2m}<y_{2m+1}$, the natural logarithm of their quotient is nonzero, so 
we are not dividing by zero.  As $y_{2m},y_{2m+1}\!\in\!(0,1)$, and as 
$\alpha\!>\!0$, both of the natural logarithms involve only positive numbers, 
so they are well defined.  It is then easily checked that if we 
first pick $b$ sufficiently large, then we can always find our required 
$c$, and thus we can construct the $g$ we desire. So we are done. \qed 

\subsection{Proving Theorem \ref{thm:3d} via the Horn-Kapranov Uniformization} 
Before outlining our proof, we will have to rapidly review $\cA$-discriminants 
and some recent bounds for the number of real roots of certain 
generalizations of binomial systems. We point out that the 
necessary background is detailed further in \cite{drrs}, so our 
key contribution here is that we extend the results of \cite{drrs} 
to higher-dimensional (reduced) discriminant varieties. 

\subsection{Background on $\cA$-Discriminants} 
First, recall that to any set $\cA\!\subset\!\Zn$, one can 
consider the family of polynomials $f$ supported on $\cA$. The 
closure of the set of coefficient values $(c_a)_{a\in \cA}$ yielding $f$ with 
a degenerate zero set in $\Csn$ is defined to be the 
$\cA$-discriminant variety $\nabla_\cA$ \cite{loeser,gkz94}. 
The $\cA$-discriminant variety indeed turns out to be an 
algebraic variety, and when $\codim \nabla_\cA\!=\!1$,  
it possesses a defining polynomial irreducible over $\Z[c_a\; | \; 
a\!\in\!\cA]$ which we will call $\Delta_\cA$. 

$\cA$-discriminant varieties turn out also to be fibered over a 
variety of even lower dimension, known as the {\bf reduced} 
$\cA$-discriminant variety, $\overline{\nabla}_\cA$. 
More to the point, the topology of $Z^*_\R(f)$ is 
constant on the {\bf real} complement of $\overline{\nabla}_\cA$.  
We now detail the necessary definitions and results. 
\begin{thm}
\label{thm:dfs}
\cite[Prop.\ 4.1]{dfs}
Given $\cA\!:=\!\{a_1,\ldots, a_m\}\!\in\!\Zn$, the discriminant locus
$\nabla_\cA$ is exactly the closure of\\
\mbox{}\hfill
$\left\{ \left[u_1t^{a_1}:\cdots:u_mt^{a_m}\right]\; \left| \; u\!:=\!(u_1,
\ldots,u_m)\!\in\!\C^m,\
\cA u \!=\!\bO, \ \sum^m_{i=1}u_i\!=\!0, \
t\!=\!(t_1,\ldots,t_n)\!\in\!\Csn\right. \right\}$. \qed\hfill\mbox{}
\end{thm}
\begin{lemma}
\label{prop:odd}
Suppose $\cA\!=\!\{a_1,\ldots,a_m\}\!\subset\!\Zn$ affinely generates $\Zn$
and $a_1\!=\!\bO$. Then there are $i_1,\ldots,i_n\!\in\!\{2,\ldots,m\}$ such
that $\det[a_{i_1},\ldots,a_{i_n}]$ is odd. \qed
\end{lemma}
                                                                                
\begin{dfn}
\label{dfn:red}
Suppose $\cA\!=\!\{a_1,\ldots,a_m\}\!\subset\!\Zn$ affinely generates $\Zn$,
has cardinality $m\!\geq\!n+2$, and $a_1\!=\!\bO$.
We call any set $C\!=\!\{i_1,\ldots,i_n\}$
with $\det[a_{i_1},\ldots,a_{i_n}]$ odd as in Lemma \ref{prop:odd} above,
an {\bf odd cell} of $\cA$. For any $n\times m$
matrix $B$, we then let $B_C$ (resp.\ $B_{C'}$) denote the submatrix of
$B$ defined by columns of $B$ with index in $C$ (resp.\ $\{2,\ldots,m\}
\setminus C$). For any vectors $v,w\!\in\!(\Cs)^m$, let us denote
their coordinate-wise product by $v\cdot w\!:=\!(v_1w_1,\ldots,v_m w_m)$.
Also let $\Gamma$ be the multivalued\footnote{The multiple values arise
from the presence of rational exponents, and the number of images of a
point is always bounded above by a constant depending only on $\cA$.}
function from $(\Cs)^m$ to $(\Cs)^{m-n-1}$ defined by
$\Gamma(y):=\frac{y_{C'}}{y_1}\cdot \left(\frac{y_C}{y_1}
\right)^{-A^{-1}_C A_{C'}}$. Finally, we define the {\bf reduced}
$\cA$-discriminant variety, $\bla\!\subset\!\C^{m-n-1}$, to be the closure
of\\
\mbox{}\hfill$\left\{ \Gamma(u) \; \left| \; u\!:=\!(u_1,
\ldots,u_m)\!\in\!(\Cs)^m,\ \cA u \!=\!\bO, \ \sum^m_{i=1}u_i\!=\!0\right.
\right\}$, \hfill\mbox{}\\
and call any connected component of $(\Rs)^{m-n-1}\setminus\bla$ a
{\bf reduced ($\cA$-)discriminant chamber}. \dia
\end{dfn}
\begin{rem}\label{rem:red}
\scalebox{.98}[1]{Since we always implicitly assume that an odd cell has
been fixed a priori for our}\linebreak
\scalebox{.9}[1]{reduced $\cA$-discriminant varieties, $\Gamma$
in fact restricts to a single-valued function from
$(\Rs)^m$ to $(\Rs)^{m-n-1}$. \dia}
\end{rem}
\begin{prop}
\label{cor:odd}
Suppose $\cA\!=\!\{a_1,\ldots,a_m\}\!\subset\!\Zn$ affinely generates $\Zn$,
has cardinality $m\!\geq\!n+2$,
and $a_1\!=\!\bO$. Also let $C$ be any odd cell of $\cA$,
let $f(x)\!:=\!\sum^m_{i=1}\delta_ix^{a_i}$ with
$\delta\!:=\!(\delta_1,\ldots,\delta_m)\!\in\!(\Rs)^m$, and let
$\bar{\delta}\!\in\!(\Rs)^m$ be
the unique vector with $\bar{\delta}_1\!=\!1$,
$\bar{\delta}_C\!=\!(1,\ldots,1)$ and
$\bar{\delta}_{C'}\!=\!\Gamma(\delta)$. Finally,
let $\bar{f}\!:=\!\sum^m_{i=1}\bar{\delta}_ix^{a_i}$ and
let $\conv \cA$ denote the convex hull of $\cA$. Then:
\begin{enumerate}
\item{ $\Gamma$ induces a surjection from the set of
connected components of\\
\mbox{}\hfill $\Pro^{m-1}_\R\setminus\left(\nabla_\cA\cup
\left\{\left.[y_1:\cdots:y_m]\!\in\!\Pro^{m-1}_\R\; \right|
\; y_1\cdots y_m\!=\!0 \right\}\right)$
\hfill\mbox{}\\
to the set of reduced $\cA$-discriminant chambers. }
\item{ If, for all facets $Q'$ of $\conv \cA$, we have that
$\#(\cA\cap Q')\!=\!n$, then $Z^*_\R(f)$ and $Z^*_\R(\bar{f})$
are diffeotopic. Furthermore, for any $f_1$ and $f_2$ with
$\bar{f}_1$ and $\bar{f}_2$ lying in the
same reduced $\cA$-discriminant chamber, $Z^*_\R(\bar{f}_1)$ and
$Z^*_\R(\bar{f}_2)$ are diffeotopic. \qed}
\end{enumerate}
\end{prop}
                                                                                
\noindent
Proposition~\ref{cor:odd} follows easily from a routine application of
the Smith normal form and the implicit function theorem. In particular, the
crucial trick is to observe that exponentiation by $A_C$, when $C$ is an odd
cell, induces an automorphism of orthants of $\Rsn$.
Our assumption on the intersection of $\cA$ with the facets of $\conv \cA$
ensures that any topological change in the
zero sets of $f$ and $\bar{f}$ (in the underlying real toric variety
corresponding to $\conv \cA$ \cite{tfulton}) occurs within $\Rsn$.

\subsection{Background on Sheared Binomial Systems} 
\begin{dfn}
\label{dfn:shear}
Suppose $\ell_1,\ldots,\ell_j\!\in\!\R[\lambda_1,\ldots,
\lambda_k]$ are polynomials of degree $\leq\!1$. We then call any system of
equations of the form
$S\!:=\!\left(1-\prod^j_{i=1} \ell^{b_{1,i}}_i(\lambda_1,\ldots,\lambda_k),
\ldots,1-\prod^j_{i=1} \ell^{b_{k,i}}_i(\lambda_1,\ldots,\lambda_k)\right)$,
with $b_{i,i'}\!\in\!\R$ for all $i,i'$, and the vectors
$(b_{1,1},\ldots,b_{1,j}),\ldots,(b_{k,1},\ldots,b_{k,j})$ linearly
independent, a {\bf $k\times k$ sheared binomial system with $j$ factors}.
We also call each $\ell_i$ a {\bf factor} of the system. A sheared binomial
system is referred to as  a
{\bf Gale Dual System} in \cite{bs}.\dia
\end{dfn}
                                                                                
\noindent
Note that our definition implies that $j\!\geq\!k$. For $j\!=\!k$,
it is easy to reduce any $k\times k$ sheared binomial system with $j$ factors
to a $k\times k$ linear system, simply by multiplying and dividing equations
(mimicking Gaussian elimination). For $j\!>\!k$, sheared polynomial systems
become much more complicated.
\begin{thm}
\cite{bs} 
\label{thm:shear} The number of non-degenerate roots
$\lambda\!\in\!\R^k$ of any $k\times k$ sheared binomial system with $n+k$
factors, and all factors positive, is bounded above by 
$(e^2+3) 2^{(k-4)(k+1)/2} n^k$, for all $k\!\geq\!1$. 
In particular, $e^2+3\!\approx\!10.38905610$. \qed
\end{thm}

\subsection{The Proof of Theorem \ref{thm:3d} } 
Let $\cT_\cA$ denote the toric variety corresponding to the convex
hull of $\cA$ \cite{tfulton}. We will separate into two cases:  
(1) those $\cA$ that are generic in the sense that 
every subset of cardinality $n+1$ has convex hull of positive 
volume, and (2) any remaining $\cA$. The most difficult case is 
Case (1), so we will restrict to this case, leaving Case (2) 
for the full version of this paper. 

By our genericity assumption, 
every facet of $\cA$ corresponds to the vertices of a simplex, and thus 
the complex zero set of any $f$ with $\supp(f)\!=\!\cA$ is always
nonsingular at infinity, relative to $\cT_\cA$ (see, e.g.,
\cite[Sec.\ 3.2]{thresh}).
By Proposition~\ref{cor:odd}, it then suffices to show that our desired
bound applies to the number of {\bf reduced} $\cA$-discriminant
chambers. Note also that by Proposition~\ref{cor:odd} and
\cite[Lemma 3.3]{drrs},
the real part of the reduced $\cA$-discriminant variety ---
$\R^3\cap\bla$ --- must be the union of a finite set of curves and the
closure of $\{\Psi(\lambda)\; | \; \lambda\!\in\!\R^2,\
\ell_1(\lambda)\cdots\ell_{n+4}(\lambda)\!\neq\!0\}$,
where $\Psi(\lambda)\!:=\!(\psi_1(\lambda),\psi_2(\lambda),\psi_3(\lambda))
\!:=\!\left(\prod^{n+4}_{i=1}\ell^{b_{1,i}}_i(\lambda),
\prod^{n+4}_{i=1}\ell^{b_{2,i}}_i(\lambda),
\prod^{n+4}_{i=1}\ell^{b_{3,i}}_i(\lambda)\right)$,
and $\ell_1,\ldots,\ell_{n+4}$ are bivariate polynomials in $\lambda$
of degree $\leq\!1$
Let $\Omega\!\subset\!\R^3$ denote the
aforementioned closure. Since curves do not disconnect connected
components of the complement of a (locally closed) real algebraic surface,
it thus suffices to focus on $\Omega$. In particular, the
connected components of $(\Rs)^3\setminus\Omega$ are (up to
the deletion of finitely many curves) exactly the reduced
$\cA$-discriminant chambers. Note also, by observing the poles of the
$\psi_i$, that $\Omega$ is the closure of the union of
no more than $(n+4)(n+5)/2$ topological disks.  
                                                                                
To count the number of connected components of $(\Rs)^3\setminus\Omega$, we
will use the classical {\bf critical points method} \cite{chigo}, combined
with our more recent tools. In particular, let us first bound the number of 
critical values of the map from $\Omega$ to $\R$ defined by $x_1$. 

A simple derivative calculation then reveals that some critical values 
are given by $\psi_3(\lambda)$, where $\lambda$ satisfies 
$\frac{\partial\psi_1}{\partial \lambda_1}
\frac{\partial\psi_2}{\partial \lambda_2}\!=\!
\frac{\partial\psi_2}{\partial \lambda_1}
\frac{\partial\psi_1}{\partial \lambda_2}$ and 
$\frac{\partial\psi_1}{\partial \lambda_1}
\frac{\partial\psi_3}{\partial \lambda_2}\!=\!
\frac{\partial\psi_3}{\partial \lambda_1}
\frac{\partial\psi_1}{\partial \lambda_2}$. 
An elementary calcluation then yields that the preceding 
$2\times 2$ system is equivalent to a $2\times 2$ polynomial 
system consisting of two polynomials of degree $\leq\!n+3$. 
So by B\'ezout's Theorem, the number of critical values is 
bounded above by $(n+3)^2$. 

Next, there are contributions from more complicated nodal singularities. 
It is then easily checked that these reduce to counting the number of 
roots $(\lambda,\lambda')\!\in\!(\Rs)^4$ of a sheared binomial system of the 
following form: $(\psi_1(\lambda),\psi_2(\lambda))\!=\!(\psi_1(\lambda),
\psi_2(\lambda))$ --- a $4\times 4$ sheared binomial system, with 
$\leq\!2n+8$ factors. By Theorem \ref{thm:shear}, and by counting 
sign conditions, we then obtain a contribution of $\frac{(2n+8)(2n+9)}{2}
\cdot (e^2+3)n^4\!=\!O(n^6)$. 

To count the number of connected components of
$(\Rs)^3\setminus\Omega$, let us now introduce planes  
$H_1,\ldots,H_N$ exactly at the locations of our preceding critical values.  
Clearly, any connected component of\\
\mbox{}\hfill$T\!:=\!(\Rs)^3\setminus(\Omega\cup H_1\cup \cdots\cup H_N)$
\hfill\mbox{}\\
is contained in a unique connected component of $(\Rs)^3\setminus\Omega$. So
it suffices to count the connected components of $T$. To do the latter,
observe that $N\!=\!O(n^6)$ and our planes thus divide $(\Rs)^3$ into 
$O(n^6)$ vertical slabs.
                                                                                
Now note that within the interior of each slab, $\Omega$
does {\bf not} intersect the $x_2$ or $x_3$ coordinate planes, and has no 
vertical tangents. So to count components of $T$ within
any particular vertical strip, we need only bound from above the number of
connected components of the complement of 
$\Omega\cup\{x_1\!=\!0\}$ within a vertical plane distinct from 
$H_1,\ldots,H_N$. This clearly reduces to the 
critical points method once more, in one dimension lower. In particular, 
with a bit of work, one is reduced asymptotically to counting the
number of real roots of a $2\times 2$ sheared binomial system 
with $n+4$ factors. So by Theorem~\ref{thm:shear} once again, and a 
sign condition count, the desired upper bound is $O(n^5)$.
Thus, each of our vertical slabs contains no more than
$O(n^5)$ connected components of $T$. Taking into account the number of
vertical strips, we thus finally arrive at our stated upper bound of 
$O(n^{11})$ for the number of connected components of $T$, so we are done.  
\qed

\section{Acknowledgements} 
We thank Martin Avenda\~{n}o for useful discussions. 

\bibliographystyle{acm}

\begin{thebibliography}{A}

\bibitem[Ave07]{avendano} Avenda\~{n}o, Martin, {\it ``The Number of 
Real Roots of a Bivariate Polynomial on a Line,''} submitted for 
publication. 

\bibitem[BV06]{basu} Basu, Saugata and Vorobjov, Nicolai N.,
{\it ``On the Number of Homotopy Types of Fibres of a Definable
Map,''} Journal of the London Mathematical Society, to appear.
Also available as Math ArXiV preprint {\tt math.AG/0605517} .

\bibitem[BS07]{bs} Bihan, Frederic and Sottile, Frank, {\it ``New
Fewnomial Upper Bounds from Gale Dual \linebreak Polynomial Systems,''}
Moscow Mathematical Journal, vol.\ 7, no.\ 3, (July--September, 2007). 

\bibitem[BRS07]{thresh} Bihan, Frederic; Rojas, J.\ Maurice; and Stella, 
Casey, {\it ``First Steps in Real Algorithmic Fewnomial Theory,''} submitted 
for publication.   

\bibitem[Can88]{pspace} Canny, John F., {\it ``Some Algebraic
and Geometric Computations in PSPACE,''} Proc.\ 20\thth ACM
Symp.\ Theory of Computing, Chicago (1988), ACM Press.

\bibitem[CG84]{chigo} Chistov, Alexander L., and Grigoriev, Dima Yu, {\it
``Complexity of Quantifier Elimination in the Theory of Algebraically
Closed Fields,''} Lect.\ Notes Comp.\ Sci.\ 176, Springer-Verlag (1984).

\bibitem[CZ02]{cohenzannier} Cohen, Paula B.\ and Zannier, Umberto, {\it
``Fewnomials and intersections of lines with real analytic subgroups in
$\mathbf{G}^n_m$,''}
Bull.\ London Math.\ Soc.\ 34 (2002), no.\ 1, pp.\ 21--32.

\bibitem[vdD98]{ominimal} van den Dries, Lou, {\it Tame topology and o-minimal
structures,} London Mathematical \linebreak Society Lecture Note Series, 248,
Cambridge University Press, Cambridge, 1998.

\bibitem[DFS05]{dfs} Dickenstein, Alicia; Feichtner, Eva Marie; 
and Sturmfels, Bernd, {\it ``Tropical Discriminants,''} 
manuscript, 2005. 

\bibitem[DRRS07]{drrs} Dickenstein, Alicia; Rojas, J.\ Maurice; 
Rusek, Korben; Shih, Justin, {\it ``Extremal Real Algebraic 
Geometry and $\cA$-Discriminants,''} Moscow Mathematical Journal, 
vol.\ 7, no.\ 3, (July--September, 2007). 

\bibitem[Ful93]{tfulton} Fulton, William, {\it
Introduction to Toric Varieties}, Annals of Mathematics Studies, no.\ 131,
Princeton University Press, Princeton, New Jersey, 1993.

\bibitem[GVZ04]{gabrielov} Gabrielov, Andrei; Vorobjov, Nikolai; and Zell, 
Thierry, {\it ``Betti Numbers of Semialgebraic and Sub-Pfaffian Sets,''} 
J.\ London Math.\ Soc.\ {\bf 69} (2004), pp.\ 27--43.  

\bibitem[GKZ94]{gkz94} Gel'fand, Israel Moseyevitch; Kapranov, Misha M.; and
Zelevinsky, Andrei V.;
{\it Discriminants, Resultants and Multidimensional Determinants,}
Birkh\"auser, Boston, 1994.

\bibitem[GLS99]{gls} Giusti, M., Lecerf, G., and Salvy, B.,
{\it ``A Gr\"obner-Free Alternative to Polynomial System Solving,''}
preprint, TERA, 1999.

\bibitem[Haa02]{haas} Haas, Bertrand, {\it ``A Simple Counter-Example
to Kushnirenko's Conjecture,''} Beitr\"age zur Algebra und Geometrie,
Vol.\ 43, No.\ 1, pp.\ 1--8 (2002).

\bibitem[Har76]{harnack} Harnack, Carl Gustav Axel, {\it ``\"Uber die
Vielfaltigkeit der ebenen algebraischen Kurven,''} Math.\ Ann.\
{\bf 10} (1876), pp.\ 189--199.

\bibitem[Kho80]{kho} Khovanski, Askold G., {\it ``On a Class of
Systems of Transcendental Equations,''} Dokl.\
Akad.\ Nauk SSSR {\bf 255} (1980), no.\ 4, pp.\ 804--807;
English transl.\ in Soviet Math.\ Dokl.\ {\bf 22} (1980),
no.\ 3.

\bibitem[Kro82]{kron} Kronecker, L., {\it ``Grundz\"uge einer
arithmetischen Theorie der algebraischen Gr\"ossen,''} J.\
reine angew.\ Math., 92:1--122, 1882.

\bibitem[Kho91]{few} \underline{\hspace{\khov}}, {\it Fewnomials,}
AMS Press, Providence, Rhode Island, 1991.

\bibitem[LRW03]{lrw} Li, Tien-Yien; Rojas, J.\ Maurice; and
Wang, Xiaoshen, {\it ``Counting Real Connected Components of Trinomial
Curves Intersections and m-nomial Hypersurfaces,''} Discrete and
Computational Geometry, 30:379--414 (2003).

\bibitem[Loe91]{loeser} Loeser, Fran\c{c}ois, {\it ``Polytopes secondaires 
et discriminants,"} S\'eminaire Bourbaki, Vol.\ 1990/91.  Ast\'erisque  
No.\ 201--203 (1991), Exp.\ No.\ 742, pp.\ 387--420 (1992). 

\bibitem[Kal03]{kaloshin} Kaloshin, V., {\it ``The existential Hilbert 16-th
problem and an estimate for cyclicity of \linebreak elementary polycycles,''}
Invent.\ Math.\ 151 (2003), no.\ 3, pp.\ 451--512.

\bibitem[OK00]{orevkov} Orevkov, S.\ Yu.\ and Kharlamov, V.\ M., {\it
``Asymptotic growth of the number of classes of real plane algebraic
curves when the degree increases,''} J.\ of Math. Sciences,  113 (2003), no.\ 5,pp.\ 666--674.

\bibitem[SL54]{descartes} Smith, David Eugene and Latham, Marcia L., {\it
The Geometry of Ren\'e Descartes,} translated from the French and Latin
(with a facsimile of Descartes' 1637 French edition),
Dover Publications Inc., New York (1954).

\end{thebibliography}

\end{document}